\documentclass[aop,MSNbibl,citesort,dvips]{arximspdf}

\doi{10.1214/11-AOP684}
\volume{40}
\issue{3}
\pubyear{2012}
\firstpage{1375}
\lastpage{1376}

\makeatletter

\newcommand{\R}{\mathbb{R}}

\makeatother

\begin{document}
\begin{frontmatter}

\title{Errata for\\
Stochastic calculus for symmetric Markov~processes}\vspace*{6pt}
\runtitle{Errata for Stochastic calculus}
\pdftitle{Errata for
Stochastic calculus for symmetric Markov processes}

\textit{Ann. Probab.} \textbf{36} (2008)
931--970

\begin{aug}
\author[A]{\fnms{Zhen-Qing} \snm{Chen}\ead
[label=e1]{zchen@math.washington.edu}},
\author[B]{\fnms{Patrick J.} \snm{Fitzsimmons}\ead
[label=e2]{pfitzsim@ucsd.edu}},
\author[C]{\fnms{Kazuhiro} \snm{Kuwae}\corref{}\ead
[label=e3]{kuwae@gpo.kumamoto-u.ac.jp}}\\ and
\author[D]{\fnms{Tu-Sheng} \snm{Zhang}\ead
[label=e4]{tusheng.zhang@manchester.ac.uk}}
\runauthor{Chen, Fitzsimmons, Kuwae and Zhang}
\affiliation{University of Washington, University of California at San
Diego,\break
Kumamoto University and University of Manchester}
\address[A]{Z.-Q. Chen\\
Department of Mathematics\\
University of Washington\\
Seattle, Washington 98195\\
USA\\
\printead{e1}}
\address[B]{P. J. Fitzsimmons\\
Department of Mathematics \\
University of California, San Diego \\
9500 Gilman Drive \\
La Jolla, California 92093-0112\\
USA\\
\printead{e2}}
\address[C]{K. Kuwae\\
Department of Mathematics and Engineering \\
Graduate School of Science and Technology \\
Kumamoto University\\
Kumamoto, 860-8555 \\
Japan\\
\printead{e3}}
\address[D]{T.-S. Zhang\\
School of Mathematics \\
University of Manchester \\
Oxford Road \\
Manchester M13 9PL \\
England, United Kingdom \\
\printead{e4}\hspace*{5.25pt}}
\end{aug}

% HISTORY:
\received{\smonth{5} \syear{2011}}

% ABSTRACT

% KEYWORDS
%

\end{frontmatter}

%s1 ###
\section{Errata}

For a correct statement of Lemma 3.1(i) in \cite{CFKZStoch}, eliminate
part (d) and replace part (c) by
\begin{longlist}[(c)]
\item[(c)] $\mathbf{P}_x(\lim_{n\to\infty}\sigma_{E\setminus G_n} <
\zeta)=0$ for $m$-a.e. $x\in E$.
\end{longlist}

To see that the original statements of parts (c) and (d) of Lemma
3.1(i) are not correct, consider the killed Brownian motion in the open
unit ball in $\R^d$, with $G_n$ the concentric open ball with radius
$1-1/n$. In Remark 4.5(iii) as well as in Lemma 4.6, for
$f\in\mathcal{F}_{\mathrm{loc}}$, the stochastic integral $\int_0^t
f(X_{s-}) \,d \Lambda(M)_s$ is in general well defined only for $t\in[0,
\zeta[$. The generalized It\^o formula in \cite{CFKZStoch}, Theorem
4.7, holds in general only for $t\in[0, \zeta[$. These corrections have
no effect on the papers \cite{CFKZPert,CFKZGenPert}, as the results
from \cite{CFKZStoch} are used in those papers only for $t\in[0,
\zeta[$.

\section*{Acknowledgments}
The authors thank Dr. Alexander Walsh for alerting us to the
problem with Lemma 3.1 in \cite{CFKZStoch}.

%suskaldyti doi

% imsref loaded by lrinkeviciute, 2011-06-23 12:54:59
%

%
\printaddresses

\end{document}